\documentclass{math_crypt_v06}

\newcommand{\M}{{\mathcal M}}

\def\QQ{{\mathbb Q}}

\def\FF{{\mathbb F}}

\def\ZZ{{\mathbb Z}}

\newcommand{\Div}{\operatorname{Div}}

\newcommand{\Jac}{\operatorname{Jac}}
\newcommand{\Prin}{\operatorname{Prin}}

\newcommand{\ord}{\operatorname{ord}}
\newcommand{\lcm}{\operatorname{lcm}}
\newcommand{\Pic}{\operatorname{Pic}}

\newcommand{\Hom}{\operatorname{Hom}}

\title{The Monodromy Pairing and \\Discrete Logarithm on the Jacobian of Finite Graphs}
\headlinetitle{DLP on the Jacobian of Finite Graphs}

\authorone{Farbod Shokrieh}
\addressone{Georgia Institute of Technology\\
Atlanta, Georgia 30332-0160}
\countryone{USA}
\emailone{shokrieh@math.gatech.edu, farbod@ece.gatech.edu}

\abstract{Every graph has a canonical finite abelian group attached to it. This group has appeared in the literature under a variety of names including the sandpile group, critical group, Jacobian group, and Picard group. The construction of this group closely mirrors the construction of the Jacobian variety of an algebraic curve. Motivated by this analogy, it was recently suggested by Norman Biggs that the critical group of a finite graph is a good candidate for doing discrete logarithm based cryptography. In this paper, we study a bilinear pairing on this group and show how to compute it. Then we use this pairing to find the discrete logarithm efficiently, thus showing that the associated cryptographic schemes are not secure. Our approach resembles the MOV attack on elliptic curves.}



\acknowledgments{I would like to thank Matthew Baker for recommending this problem to me, and for many helpful discussions. Thanks also to Xander Faber, Dino Lorenzini, and Prasad Tetali for their helpful comments.\\ The author's work was supported in part by NSF grants DMS-0901487 and CCF-0902717, during the writing of this article.}

\firstpage{43}                                %
\doi{10.1515}                                    %
\communicated{Simon R. Blackburn}                           %
\received{30 July, 2009}                               %
\revised{2 March, 2010}                                %

\begin{document}

\section{Introduction}
\label{IntroSection}

\subsection{Overview}
\label{OverviewSection}

 Every graph has a canonical finite abelian group whose order is the number of spanning trees of the graph. This group has appeared in the literature under many different names; in theoretical physics it was first introduced as the ``abelian sandpile group'' or ``abelian avalanche group''  in the context of self-organized critical phenomena (\cite{BTW88,Dhar90,Gabrielov93}). In arithmetic geometry, this group appeared as the ``group of components'' in the study of degenerating algebraic curves (\cite{Lorenzini89}). In algebraic graph theory this group appeared under the name ``Jacobian group'' or ``Picard group'' in the study of flows and cuts in graphs (\cite{BacherHN97}). The study of a certain chip-firing game on graphs led to the definition of this group under the name ``critical group'' (\cite{Biggs97,Biggs99}).

The construction of this group closely mirrors the construction of the Jacobian variety of an algebraic curve. Motivated by this analogy, Norman Biggs in \cite{Biggs07} suggests that the Jacobian of a finite graph (which he calls the ``critical group'') might be suitable for discrete logarithm based cryptography.

In this paper, we study the discrete logarithm problem on the Jacobian of finite graphs. Our main result is an algorithm to efficiently compute discrete logarithms on these groups. Therefore, unlike elliptic curves and Jacobian varieties, one can not use the Jacobian of finite graphs for cryptographic purposes. It is an intriguing problem whether the fact that discrete logarithm can be done efficiently might have any algorithmic applications. Our algorithm uses a bilinear pairing, which we call the monodromy pairing, on this group. This approach is similar to the MOV attack on elliptic curves. For our application, we study the monodromy pairing and show how to compute it.

\subsection{Related work}
 The order of the Jacobian group is the number of spanning trees of the graph (\cite{Biggs97}). Hence, the order of the group can be computed by the famous Matrix-Tree formula of Kirchhoff.

Finite graphs and algebraic curves behave similarly in many respects. Recently, there have been an increasing number of papers pursuing this analogy. Some relationship between elliptic curves and chip-firing games on graphs is noticed in \cite{Musiker07}. In \cite{BN1,BN2} a version of the famous Riemann-Roch theorem is proved for finite graphs, a discrete analogue of holomorphic maps between Riemann surfaces is introduced, and a graph-theoretic Riemann-Hurwitz formula is proved. A Torelli's theorem for graphs is proved in \cite{Artamkin06,CaporasoViviani09}.
The relationship between graph theory and algebraic geometry goes beyond a simple analogy. For example, Mikhalkin and Zharkov in \cite{MK08} prove that an (abstract) ``tropical curve'' is simply a connected ``metric graph''.

\subsection{Previous work.} Norman Biggs in \cite{Biggs07} constructs a family of graphs with cyclic Jacobian groups, to be potentially used for cryptography. The problem of finding families of graphs with cyclic Jacobian groups is subsequently studied in \cite{Lorenzini08,ChenYe09,Musiker07}. These provide examples of cyclic Jacobian groups with appropriate order, so that discrete logarithm problem cannot be solved by the known purely group-theoretic methods.

In \cite{Blackburn08}  Blackburn addresses the discrete logarithm problem for the particular family of graphs constructed by Biggs in \cite{Biggs07}. It is fairly clear that methods presented in \cite{Blackburn08}, with some minor modifications, can also be applied to the general case. Our method is quite different from Blackburn's method, and our algorithm in \S\ref{DLPSection} works for any graph.

To our knowledge, the monodromy pairing was first introduced by Bosch and Lorenzini in \cite{BoschLorenzini02}. We have not found an easy-to-compute formula, like (\ref{MonodromyNonCanonicalDef}), in the literature.

\medskip

\medskip

The paper proceeds as follows. In \S\ref{DefinitionsSection} we provide the relevant definitions. The monodromy pairing is studied in \S\ref{BilinearSection}. Using the monodromy pairing, we give our discrete logarithm algorithm in \S\ref{DLPSection}. Further remarks and results are outlined in \S\ref{ComplementsSection}. Appendix~\ref{PairingProofsSection} contains a new proof of Theorem~\ref{PairingTheorem}.

\section{Definitions}
\label{DefinitionsSection}
\subsection{Notation and Terminology}
\label{NotationSection}
Throughout this paper, a {\em graph} means a finite, unweighted multigraph with no loops. All graphs are assumed to be connected. For a graph $G$, the set of vertices is denoted by $V(G)$, and the set of edges is denoted by $E(G)$. Throughout this paper, $n$ and $m$ denote the number of vertices and edges, respectively.

Let $\{ v_1,\ldots, v_n \}$ be an ordering of $V(G)$. With respect to this ordering, the {\em Laplacian matrix} $Q$ associated to $G$ is the $n \times n$
matrix $Q =(q_{ij})$, where $q_{ii}$ is the degree of vertex $v_i$, and $-q_{ij}$ ($i \neq j$) is the number of edges connecting $v_i$ and $v_j$.
It is well-known (and easy to verify) that $Q$ is symmetric, has rank $n-1$, and the kernel of $Q$ is spanned by $\mathbf{1}$, the all-one vector\footnote{Remember that $G$ has no loops.} (see, e.g., \cite{BiggsBook93,Bollobas98}).

\subsection{The Jacobian of a finite graph}
 \label{JacobianDefinitionSubsection}

Let $\Div(G)$ be the free abelian group generated by $V(G)$. One can think of elements of $\Div(G)$ as formal integer linear combination of vertices
\[
\Div(G)=\{\sum_{v \in V(G)} a_v (v) \; : \;  a_v  \in \ZZ \} \ .
\]
By analogy with the algebraic curve case, elements of $\Div(G)$ are called {\em divisors} on $G$. For a divisor $D$, the coefficient $a_v$ of $(v)$ in $D$ is denoted by $D(v)$.

We define by $\M(G) = \Hom(V(G), \ZZ) $ the abelian group consisting of all integer-valued functions on the vertices. One can think of $\M(G)$ as analogous to the group $\M(X)^{\times}$ of nonzero meromorphic functions on an algebraic curve $X$.

For $f \in \M(G)$, $\operatorname{div}(f) \in \Div(G)$ is given by the formula
\[
\operatorname{div}(f) = \sum_{v \in V(G)} {\ord_v(f)} (v) \ ,
\]
where
\[
\ord_v(f)=\sum_{\{v,w\}  \in E(G)} (f(v) - f(w)) \ .
\]
Consider the group homomorphism $\deg : \Div(G) \to \ZZ$ defined by $\deg(D) = \sum_{v \in V(G)} D(v)$. Denote by $\Div^0(G)$ the kernel of this homomorphism, consisting of {\em divisors of degree zero}. Define $\Prin(G)=\{ \operatorname{div}(f) \in \Div(G) \; : \;  f \in \M(G)\}$ to be the group of {\em principal divisors}.
\begin{lemma}
$\Prin(G) \subseteq \Div^0(G)$, and both $\Prin(G)$ and $\Div^0(G)$ are free $\ZZ$-modules of rank $n-1$.
\end{lemma}
A proof is given in \cite{Biggs97}. As a corollary, the quotient group
\[
\Jac(G) = \Div^0(G) / \Prin(G)
\]
is well-defined and is a finite abelian group. Following \cite{BacherHN97}, it is called the {\em Jacobian} or the {\em Picard}\footnote{Another appropriate notation is $\Pic^0(G)$.} group of $G$.

The following lemma is a direct consequence of Kirchhoff's famous {\em Matrix-Tree Theorem} \cite{Kirchhoff1847} (see also \cite{BacherHN97,Biggs97}).
\begin{lemma} \label{JacobianLemma}
 The order of the group $\Jac(G)$ is equal to the number of {\em spanning trees} in $G$, which we denote by $\kappa(G)$.
\end{lemma}

Following \cite{BN1}, for $D_1, D_2 \in \Div(G)$, we say that $D_1$ is equivalent to $D_2$, and write $D_1 \sim D_2$, if $D_1-D_2$ is a principal divisor.


\section{A bilinear pairing on the Jacobian of finite graphs}
\label{BilinearSection}

\subsection{Generalized inverses}
 \label{GeneralizedInverses}
A matrix can have an inverse only if it is square and its columns (or rows) are linearly independent. But one can still get ``partial inverse'' of any matrix.
\begin{definition}
Let $A$ be a matrix (not necessarily square). Any matrix $L$ satisfying $ALA=A$ is called a generalized inverse of $A$.
\end{definition}
It is somehow surprising that for every matrix $A$ there exists at least one generalized inverse. In fact, more is true; every matrix has a unique {\em Moore-Penrose pseudoinverse}\footnote{The Moore-Penrose pseudoinverse of $A$ is a generalized inverse of $A$ with three extra properties; see \cite{BenIsrael03} for an extensive study of the subject.}.

\medskip

Let $Q$ be the Laplacian matrix of a connected graph. Since its rank is $n-1$, it cannot have an inverse. But there are many ways to obtain generalized inverses:

\begin{example} \label{GeneralizeInvExample1}
Fix an integer $1 \leq i \leq n$. Let $Q_i$ be the $(n-1) \times (n-1)$ matrix obtained from $Q$ by deleting $i^{\rm th}$ row and $i^{\rm th}$ column. Then $Q_i$ is a full rank matrix and has an inverse $Q_i^{-1}$. Let $L_{(i)}$ be the $n \times n$ matrix obtained from $Q_i^{-1}$ by inserting a zero row after the $(i-1)^{\rm th}$ row and inserting a zero column after the $(i-1)^{\rm th}$ column. Then $L_{(i)}$ is a generalized inverse of $Q$. One can check that
\[QL_{(i)}=I+R_{(i)} \ , \]
where $I$ is the identity matrix, and $R_{(i)}$ has $-1$ entries in the $i^{\rm th}$ row and is zero everywhere else. As $R_{(i)}Q=0$, we get $QL_{(i)}Q=Q$.
\end{example}

\begin{example}\label{GeneralizeInvExample2}
Let $J$ be the $n \times n$ all-one matrix. Then $Q+\frac{1}{n}J$ is nonsingular and $Q^+= (Q+\frac{1}{n}J)^{-1}-\frac{1}{n}J$ is a generalized inverse of $Q$. In fact it is the unique {\em Moore-Penrose pseudoinverse} of $Q$; it is easy to check $QQ^{+}=Q^{+}Q=I-\frac{1}{n}J$ and $Q^{+}QQ^{+}=Q^{+}$.
\end{example}
These examples show that computing a generalized inverse $L$ takes time at most $O(n^{\omega})$, where $\omega$ is the exponent for matrix multiplication.

\subsection{The monodromy pairing}
\label{MonodromySubsection}
A kind of graph-theoretic analogue of Weil pairing on the (principally polarized) Jacobian of an algebraic curve is provided by a certain bilinear pairing on $\Jac(G)$, which we define in this section\footnote{The monodromy pairing is symmetric, while the Weil pairing is skew-symmetric.}.

For $D_1, D_2$ in $\Div^0(G)$, let $m_1$ and $m_2$ be integers such that $m_1D_1=\operatorname{div}(f_1)$ and $m_2 D_2=\operatorname{div}(f_2)$ are principal; these exist because $\Jac(G)$ is a finite group. One can easily show that
\begin{equation}\label{DivPairSymmetry}
\frac{1}{m_2}\sum_{v\in V(G)}{D_1(v) f_2(v)}=\frac{1}{m_1}\sum_{v\in V(G)}{f_1(v) D_2(v)} \ .
\end{equation}
The pairing $\langle  \cdot \, , \cdot \rangle: \; \Div^0(G) \times \Div^0(G) \rightarrow \QQ$ defined by
\begin{equation}
\langle D_1 , D_2 \rangle = \frac{1}{m_2}\sum_{v\in V(G)}{D_1(v) f_2(v)}
\end{equation}
is symmetric and bilinear. This pairing descends to a well-defined pairing on $\Jac(G)$. We use the notation $\overline{D}$ for an element of $\Jac(G)$, if $D$ is a lift of that element in $\Div^0(G)$.

\begin{theorem}\label{PairingTheorem}
The pairing $\langle  \cdot \, , \cdot \rangle : \; \Jac(G) \times \Jac(G) \rightarrow \QQ/\ZZ$ defined by
\begin{equation}\label{MonodromyCanonicalDef}
\langle \overline{D}_1 , \overline{D}_2 \rangle = \frac{1}{m_2}\sum_{v\in V(G)}{D_1(v) f_2(v)} \pmod{\ZZ} \ ,
\end{equation}

where $m_2D_2=\operatorname{div}(f_2)$, is a well-defined, symmetric, bilinear, non-degenerate pairing on $\Jac(G)$.
\end{theorem}

\medskip

This theorem, in a slightly different language, is proved in \cite{BoschLorenzini02}. We give a more elementary proof in Appendix~\ref{PairingProofsSection}.
\begin{definition}
 We call the pairing described in Theorem~\ref{PairingTheorem} the {\em monodromy pairing} (see remark 1 in \S\ref{ComplementsSection} for this terminology).
\end{definition}

\begin{remark}
Let $\Phi$ be a finitely generated abelian group. A symmetric bilinear pairing $\langle \, \cdot \, , \cdot \rangle : \; \Phi \times \Phi \rightarrow \QQ/\ZZ$ is called {\em non-degenerate} (or {\em regular}) if the group homomorphism $\Phi \rightarrow \Hom_{\ZZ}(\Phi, \QQ/\ZZ)$ defined by $x \mapsto  \langle x, \cdot \rangle$ is injective. If it is an isomorphism, it is called {\em perfect} (or {\em unimodular}). If a pairing on a finitely generated abelian group is non-degenerate, then it is automatically  perfect\footnote{Moreover, the group is torsion in this situation.} (see \cite{Durfee77}). For a finite abelian group $\Phi$, this fact is immediate; there exists a (non-canonical) isomorphism between $\Phi$ and its {\em Pontryagin dual} $\Hom_{\ZZ}(\Phi, \QQ/\ZZ)$ (see, e.g., page 167 of \cite{DummitFoote04}).
\end{remark}

Let $\{ v_1,\ldots, v_n \}$ be an ordering of $V(G)$. Let $Q$ be the Laplacian matrix with respect to this ordering. This ordering gives an isomorphism between abelian groups $\Div(G)$, $\M(G)$, and the $\ZZ$-module of $n \times 1$ column vectors having integer coordinates. Under these isomorphisms the operator $\operatorname{div} : \M(G) \to \Div(G)$ coincides with the $\ZZ$-module homomorphism $Q: \ZZ^n \to \ZZ^n$. More specifically, if $[D]$ denotes the column vector corresponding to $D \in \Div(G)$, and  $[f]$ denotes the column vector corresponding to $f \in \M(G)$, then $[\operatorname{div}(f)] = Q [f]$.

The given definition of the monodromy pairing in (\ref{MonodromyCanonicalDef}) is canonical. However, the following proposition simplifies the proof of Theorem~\ref{PairingTheorem}. Moreover, it shows how one can compute the monodromy pairing in practice.

\begin{proposition}\label{PairingProposition}
Let $L$ be {\em any} generalized inverse of the Laplacian matrix $Q$. Then the monodromy pairing is given by
\begin{equation}\label{MonodromyNonCanonicalDef}
\langle \overline{D}_1 , \overline{D}_2 \rangle = [D_1]^{T} L [D_2] \pmod{\ZZ} \ .
\end{equation}
\end{proposition}
\begin{proof}
By definition $m_i[D_i]=[\operatorname{div}(f_i)]=Q[f_i]$ for $i=1,2$. The result follows from the following computations. All equalities are mod $\ZZ$
\[
\begin{aligned}
\langle \overline{D}_1 , \overline{D}_2 \rangle &= \frac{1}{m_2}\sum_{v\in V(G)}{D_1(v) f_2(v)} \\
   &= \frac{1}{m_2} [D_1]^{T}[f_2]  \\
   &= \frac{1}{m_1m_2}(Q[f_1])^{T} [f_2] \\
   &= \frac{1}{m_1m_2}[f_1]^{T}Q[f_2]  \\
   &= \frac{1}{m_1m_2}[f_1]^{T}QLQ[f_2]  \\
   &= \frac{1}{m_1m_2}(Q[f_1])^{T}L(Q[f_2]) \\
   &= [D_1]^{T} L [D_2]  \pmod{\ZZ} \ . \\
\end{aligned}
\]
\end{proof}

We emphasize that {\em any} generalized inverse of the Laplacian matrix can be used in (\ref{MonodromyNonCanonicalDef}).


\section{Discrete Logarithm Problem on the Jacobian of a finite graph}
\label{DLPSection}
Let $(\Phi,+)$ be a cyclic group. The {\em Discrete Logarithm Problem} (DLP) can be stated as:
\begin{center}
Given $g,h \in \Phi$ with $x \cdot g=h$ for some integer $x$, compute $x \mod{\ord(g)}$.
\end{center}

\medskip

In this section we use the monodromy pairing to solve the DLP for the Jacobian of a finite graph $G$ when $\Jac(G)$ is cyclic.

In our context, we assume the elements of $\Jac(G)$ are presented by some (arbitrary) lifts  in $\Div^0(G)$. Also, we assume\footnote{There are several efficient methods to find a generator; we omit the details here.} a generator $\mathbf{g}$ of the cyclic group $\Jac(G)$ is known. We can compute and save a generalized inverse $L$ of $Q$ as outlined in \S\ref{GeneralizedInverses}.

\medskip

\begin{algorithm}
(\textbf{DLP on $\Jac(G)$})
\newline \textbf{Input}: $D,D' \in \Div^0(G)$ such that $\overline{D'}=x \cdot \overline{D}$ in $\Jac(G)$
\newline \textbf{Output}: $x \mod{\ord(\overline{D})}$, the order of $\overline{D} \in \Jac(G)$.
\begin{itemize}
\item[(1)] Compute $\langle \overline{D} , \mathbf{g} \rangle=r+\ZZ$ and $\langle \overline{D'} , \mathbf{g} \rangle=r'+\ZZ$ using formula (\ref{MonodromyNonCanonicalDef})
\item[(2)] Solve the Diophantine equation $r'=rx+y$ (for variables $x,y \in \ZZ$) by clearing the denominators of $r$ and $r'$ and using the extended Euclidean algorithm, to get $x \mod{\ord(\overline{D})}$.
\end{itemize}
\end{algorithm}

\medskip

\textbf{Analysis of the algorithm.}
Since the monodromy pairing is bilinear, we have $\langle \overline{D'} , \mathbf{g} \rangle=x \langle \overline{D} , \mathbf{g} \rangle$ , or $r'=rx$ in $\QQ/\ZZ$. We still need to prove that solving the Diophantine equation precisely gives $x$ modulo the order of $\overline{D}$ in $\Jac(G)$.

\medskip

\begin{lemma}\label{OrderLemma}
Let  $\mathbf{g}$ be a generator of the cyclic group $\Jac(G)$. Let $\mathbf{h}$ be any element of  $\Jac(G)$. If $\langle \mathbf{h} , \mathbf{g} \rangle= \frac{a}{b}+\ZZ$ ($a,b \in \ZZ$, $\gcd(a,b)=1$) then $b$ is precisely the order of $\mathbf{h}$ in  $\Jac(G)$.
\end{lemma}
\begin{proof}
Let $\gamma$ be the order of $\mathbf{h}$ in  $\Jac(G)$. By bilinearity of the monodromy pairing, $\frac{\gamma a}{b}+\ZZ=\gamma \langle \mathbf{h} , \mathbf{g} \rangle=\langle \gamma
\cdot \mathbf{h} , \mathbf{g} \rangle=\langle \mathbf{0} , \mathbf{g} \rangle=\ZZ$, and therefore $b| \gamma$.
\newline On the other hand, $\langle b \cdot \mathbf{h} , \mathbf{g} \rangle=b \langle \mathbf{h} , \mathbf{g} \rangle=a+\ZZ=\ZZ$. Since $\Jac(G)$ is cyclic and the monodromy pairing is bilinear, all elements of $\Jac(G)$ must pair trivially with $b \cdot \mathbf{h}$. By non-degeneracy of the monodromy pairing we get $b \cdot \mathbf{h}=\mathbf{0}$, which means $\gamma | b$. Therefore $\gamma=b$.
\end{proof}

\medskip

Now we can show that the algorithm precisely gives $x$ mod the order of $\overline{D}$ in $\Jac(G)$; since $\overline{D'}=x\overline{D}$, order of $\overline{D'}$ divides the order of $\overline{D}$. By Lemma~\ref{OrderLemma}, for $r=\frac{a}{b}$ ($a,b \in \ZZ$, $\gcd(a,b)=1$), b is the order $\overline{D}$ in $\Jac(G)$. Multiplying by $b$ clears the denominators  in $r'=rx+y$, and we get $ax+by=c$, for some integer $c$. It is an elementary fact that the linear Diophantine equation $ax+by=c$ (with $\gcd(a,b)=1$) has solution, and $x$ is determined mod $b$.

\medskip

\textbf{Running time.}
The monodromy pairings in step (1) can be computed using the formula (\ref{MonodromyNonCanonicalDef}) given in Proposition~\ref{PairingProposition}. For this, we need to have a generalized inverse $L$; this takes time at most $O(n^{\omega})$, where $\omega$ is the exponent for matrix multiplication. Notice that this computation is done only once. Each monodromy computation in step (1) can be done using $O(n^2)$ operations (multiplication and addition). Step (2) can also be done in $O(n^2)$ operations.

For bit complexity, note that if $L=L_{(i)}$ as in Example~\ref{GeneralizeInvExample1}, one can easily see that the denominators appearing in the generalized inverse are annihilated by the exponent of the Jacobian group. The exponent is bounded above by the number of spanning trees of the graph. If we allow at most $c$ parallel edges then there are at most $c^{n-1}\cdot n^{n-2}$ spanning trees. Therefore, the denominators can be represented by $O(n \cdot \log{cn})$ bits. Moreover, one can also show that the absolute value of the entries of $L_{(i)}$ are bounded above by $R_{\max}$, the maximum effective resistance between any two vertices of the graph\footnote{$R_{\max}$ is always bounded above by the diameter of the graph, but often is much smaller than the diameter.}. Therefore all integers in the algorithm can be represented in $O(n \cdot \log{cn})$ bits.

\begin{remark}
\begin{itemize}
\item[(i)] The number of spanning trees of a graph on $n$ vertices can be exponential (or more) in $n$. For example, Biggs in \cite{Biggs07} constructs a family of graphs with cyclic Jacobian groups of order exponential in $n$. Therefore, our algorithm truly beats the known group-theoretic methods of solving DLP on finite groups.
\item[(ii)] It follows from the above discussion that, for any finite cyclic group $\Phi$, whenever one can construct an efficiently computable perfect bilinear pairing $\langle \, \cdot \, , \cdot \rangle : \; \Phi \times \Phi \rightarrow \QQ/\ZZ$, then the DLP is easy on $\Phi$. 
\end{itemize}
\end{remark}

\begin{remark}[DLP on the Critical group of finite graphs]\label{CriticalGroupRemark}
Fix a vertex $q$. In \cite{Biggs07} each element of the Jacobian group is presented by a canonical (relative to the base vertex $q$) lift in $\Div^0(G)$, which is called the {\em critical configuration} (based at $q$), and is defined by a certain chip-firing game on the graph. It is known that in each equivalence class of divisors there is a unique such critical configuration. {\em $q$-reduced divisors} (or {\em $G$-parking functions} based at $q$) provide another set of canonical elements for equivalence classes (see, e.g., \cite{FarbodReduced09} and references therein). Hence, the group law on the $\Jac(G) = \Div^0(G) / \Prin(G)$ induces a group law on the set of $q$-critical configurations, or the set of $q$-reduced divisors\footnote{In particular, cardinality of these sets are equal to the number of spanning trees.}. Biggs (\cite{Biggs99}) calls the former set with the induced group law the {\em critical group} $\operatorname{K}(G)$ of the graph, and suggests in \cite{Biggs07} that the DLP is hard on the critical group. We note that the algorithm given in this paper works for {\em any} lift of elements of $\Jac(G)$ to $\Div^0(G)$, and therefore it also solves the DLP on the critical group, as well as the ``reduced divisors group''. Some related algorithmic questions are studied in \cite{FarbodReduced09} and \cite{Heuvel01}.
\end{remark}


\section{Concluding remarks}
\label{ComplementsSection}
We conclude with some remarks.
 \begin{itemize}
\item[1.] The pairing described in Theorem~\ref{PairingTheorem} is called the {\em monodromy} or {\em Grothendieck's} pairing for the following reason. If $K$ is the field of fractions of a strictly Henselian discrete valuation ring $R$, then a theorem of Raynaud \cite{RaynaudPicard} asserts that the component group $\Phi_J$ of the N\'{e}ron model of the Jacobian $J$ of a semistable curve $X/K$ is isomorphic to the Jacobian of the dual graph $G$ of the special fiber of any semistable regular model for $X$ over $R$. Under the isomorphism provided by Raynaud's theorem, the pairing on $\Jac(G)$ which we described in \S\ref{MonodromySubsection} corresponds to Grothendieck's monodromy pairing on $\Phi_J$ (see \cite{BoschLorenzini02}).
\item[2.] By Abel's theorem for graphs (see \cite{BacherHN97}), there is a canonical isomorphism \[\Div^0(G) / \Prin(G) \cong H_1(G,\ZZ)^{\#}/H_1(G,\ZZ)\] where  $H_1(G,\ZZ)^{\#}$ denotes the dual of the cycle lattice $ H_1(G,\ZZ)$ with respect to the standard inner product on the 1-chain group $ C_1(G,\ZZ)$. It can be shown that under this canonical isomorphism, the monodromy pairing on $\Div^0(G) / \Prin(G)$ corresponds to the negative of the {\em discriminant form}\footnote{If $\Lambda$ is an integral lattice (i.e., a free $\ZZ$-module of finite rank endowed with a non-degenerate $\ZZ$-valued symmetric bilinear form), then the dual lattice $\Lambda^{\#}$ contains $\Lambda$ as a finite index subgroup, and the quotient group $\Lambda^{\#}/\Lambda$ (called the {\em discriminant group} of the lattice) inherits in a natural way a non-degenerate
$\QQ/\ZZ$-valued symmetric bilinear form, called the discriminant form (see \cite{Nikulin79} for more details).} on $H_1(G,\ZZ)^{\#}/H_1(G,\ZZ)$. This and some relevant results will appear in a subsequent paper by the author.
\item[3.] Our approach to solve the DLP on the Jacobian of finite graphs resembles the MOV attack of Menezes, Okamoto, and Vanstone \cite{MOV93} for the DLP on elliptic curves. However, because the target group of the monodromy pairing is $\QQ/\ZZ$ (instead\footnote{We also note that the monodromy pairing is symmetric, while the Weil pairing is skew-symmetric.} of $\FF_{q^{\alpha}}$), and because of Lemma~\ref{OrderLemma}, we get a deterministic polynomial-time solution for cyclic Jacobian (instead of a probabilistic polynomial-time reduction to the DLP in the group $\FF_{q^{\alpha}}^{\times}$).
\item[4.] If $\Jac(G)$ is not cyclic, then one can still use the monodromy pairing, and solve the DLP efficiently. Given a set $\{\mathbf{g}_1 , \cdots , \mathbf{g}_s\}$ generating $\Jac(G)$, the idea is to compute $\langle \overline{D} , \mathbf{g}_i \rangle=r_i+\ZZ$ and $\langle \overline{D'} , \mathbf{g}_i \rangle=r'_i+\ZZ$ for $1 \leq i \leq s$, and then solve all Diophantine equations $r'_i=r_ix+y_i$ (for variables $x,y_i \in \ZZ$) by clearing the denominators of $r_i$ and $r'_i$ and using the extended Euclidean algorithm, to get $x \mod{b_i}$ (where $r_i=a_i/b_i$ , $\gcd(a_i,b_i)=1$ ). Eventually, using Chinese Remainder Theorem, one can compute $x \mod{(\ord(\overline{D}) = \lcm{(b_1 \cdots, b_s)})}$. Interested reader can work out the details.
\item[5.] We have found at least two other methods of solving the DLP in this context. One method is essentially applying the independent work of Blackburn (\cite{Blackburn08}) to arbitrary graphs. 

\item[6.] It is worth investigating how the given solution to the DLP for the Jacobian of finite graphs can relate to the DLP for the Jacobian of algebraic curves.
\item[7.] It is an intriguing problem whether the fact that discrete logarithm can be done efficiently might have any algorithmic applications. Also the fact that the Jacobian group is actually a bilinear group with an efficiently computable pairing might have other algorithmic applications.
\end{itemize}

\appendix
\section{Proof of Theorem~\ref{PairingTheorem}}
\label{PairingProofsSection}

Here we outline an elementary proof of Theorem~\ref{PairingTheorem}. We choose an ordering of $V(G)$ and use the formula (\ref{MonodromyNonCanonicalDef}).
\newline {\em Pairing is bilinear}. This is obvious!
\newline {\em  Pairing is symmetric}. This follows from (\ref{DivPairSymmetry}). Alternatively, if $L$ is any generalized inverse of $Q$ then $L^{T}$ is also a generalized inverse of $Q$ (because $Q$ is symmetric) and we have
\[
\begin{aligned}
\langle \overline{D}_2 , \overline{D}_1 \rangle &= [D_2]^{T} L [D_1]  \pmod{\ZZ} \\
&= ([D_1]^{T} L^{T} [D_2])^{T} \pmod{\ZZ} \\
&= [D_1]^{T} L^{T} [D_2]  \pmod{\ZZ} \\
&= \langle \overline{D}_1 , \overline{D}_2 \rangle \ . \\
\end{aligned}
\]
\newline {\em  Pairing is well-defined}. Let $D_2$ and $D'_2$ in $\Div^0(G)$ be two different lifts of $\overline{D}_2 \in \Jac(G)$. Then they differ by a principal divisor $[D'_2]=[D_2]+Q[g]$ for some $g \in \M(G)$. Let $m_1[D_1]=Q[f_1]$. Then
\[
\begin{aligned}
\ [D_1]^{T} L [D'_2] &= [D_1]^{T} L [D_2]+ [D_1]^{T} L Q[g] \\
&=  [D_1]^{T} L [D_2]+\frac{1}{m_1}[f_1]^{T}QLQ[g]\\
&=  [D_1]^{T} L [D_2]+\frac{1}{m_1}[f_1]^{T}Q[g]\\
&=  [D_1]^{T} L [D_2]+[D_1]^{T}[g] \ .\\
\end{aligned}
\]
So $[D_1]^{T} L [D'_2]=[D_1]^{T} L [D_2] \pmod{\ZZ}$. By symmetry the same is true for different lifts of $D_1$.
\newline {\em  Pairing is non-degenerate}. We should show that if $D_1 \in \Div^0(G)$ be such that \[\langle \overline{D}_1, \cdot \rangle \equiv 0 \pmod{\ZZ}\] then $D_1$ is a principal divisor. Let $\mathbf{x}=L^{T}[D_1] \in \QQ^n$. If $\langle \overline{D}, \cdot \rangle \equiv 0 \pmod{\ZZ}$, then $\mathbf{x}^{T} \mathbf{u} \in \ZZ$, for any zero-sum column vector $\mathbf{u} \in \ZZ^n$. Substituting $\mathbf{e}_i-\mathbf{e}_1$ for $\mathbf{u}$ (where $\mathbf{e}_i$ denotes the vector with a $1$ in the $i^{\rm th}$ coordinate and $0$'s elsewhere), we get $\mathbf{x}=r\mathbf{1}+\mathbf{v}$ for some $r \in \QQ$ and $\mathbf{v}\in \ZZ^n$ ($\mathbf{1}$ denotes the all-one vector). Multiplying by $Q$, we have $Q\mathbf{x}=rQ\mathbf{1}+Q\mathbf{v}=Q\mathbf{v}$ or $QL^{T}[D_1]=Q\mathbf{v}$. Using $m_1[D_1]=Q[f_1]$ and the fact that $L^{T}$ is also a generalized inverse of $Q$, we get
\[
\begin{aligned}
QL^{T}[D_1]&=\frac{1}{m_1}QL^{T}Q[f_1]\\
&=\frac{1}{m_1}Q[f_1]\\
&=[D_1] \ .\\
\end{aligned}
\]
 Therefore we have shown  $[D_1]=Q\mathbf{v}$ for some $\mathbf{v}\in \ZZ^n$, which means $D_1$ must be a principal divisor.
\qed


\bibliography{jmc09}

\end{document}